\documentclass{article}

\usepackage{graphicx}
\usepackage[T1]{fontenc}
\usepackage[latin1]{inputenc}
\usepackage[french]{babel}
\usepackage{array}

\usepackage{amsmath,amsfonts,amssymb}
\usepackage{amsfonts}
\usepackage{stmaryrd}

\let\ds=\displaystyle

\let\equivalent=\Longleftrightarrow
\let\to=\rightarrow
\newcommand{\be}{\begin{enumerate}}
\newcommand{\ee}{\end{enumerate}}
\newcommand{\bi}{\begin{itemize}}
\newcommand{\ei}{\end{itemize}}

\def \abs#1{{\left|#1\right|}}
\def \pa#1{{\left(#1\right)}}

\catcode`\é=\active
\defé{\'e}

\catcode`\è=\active
\defè{\`e}

\catcode`\ê=\active
\defê{\^e}

\catcode`\û=\active
\defû{\^u}

\catcode`\â=\active
\defâ{\^a}

\catcode`\à=\active
\defà{\`a}

\catcode`\ù=\active
\defù{\`u}

\catcode`\ç=\active
\defç{\c c}

\catcode`\î=\active
\defî{\^{\i}}

\catcode`\ï=\active
\defï{\"{\i}}

\catcode`\ô=\active
\defô{\^o}

\catcode`\é=\active
\defé{\'e}

\def\bs{\bigskip}
\def\ms{\medskip}

\font\ineg=msam8
\def\ie{\mathrel{\hbox{\ineg 6}}} 
\def\se{\mathrel{\hbox{\ineg >}}} 

\font\matcinq=msbm5
\font\matsept=msbm7
\font\matdix=msbm10
\newfam\matfam \scriptscriptfont\matfam=\matcinq \scriptfont\matfam=\matsept \textfont\matfam=\matdix
\def\mat{\fam\matfam}
\def\N{{\mat N}} \def\Z{{\mat Z}}
 
\def\R{{\mat R}} \def\C{{\mat C}}

\newcommand{\bib}[2]{\hbox{\hbox to 12mm{[#1] :\hfill} \hfill \hbox to 144mm{\vtop{\hsize=144mm#2\vfill\ms}\hfill}}}

\let\phi=\varphi

\def\sin{{\rm sin}}
\def\cos{{\rm cos}}
\def\sh{{\rm sh}}
\def\ch{{\rm ch}}

\def\twp{\tilde{\wp}} 
\def\twpa{\varcurlyvee} 

\begin{document}

\begin{center}
\LARGE Elliptic functions revisited
\end{center}
\ms

\begin{center}
\large Jean-Christophe Feauveau
\footnote{Jean-Christophe Feauveau,\\
Professeur en classes préparatoires au lycée Bellevue,\\
135, route de Narbonne BP. 44370, 31031 Toulouse Cedex 4, France,\\
email : Jean-Christophe.Feauveau@ac-toulouse.fr
}

\end{center}

\begin{center}
\large January 30th, 2017
\end{center}

\bs

\textsc{Abstract.}

Elliptic functions are largely studied and standardized mathematical objects. The two usual approaches are due to Jacobi and Weierstrass.
\ms

From a contour integral which allowed us to unify many summation formulae (Euler-MacLaurin, Poisson, Voronoï or Circle formulae), we will find the entirety of the elliptic functions, proposed either in the shape of Jacobi or Weierstrass. But with one translation which appears in their natural form.
\ms

What could seem a defect will lead us to a renormalisation of the elliptic functions making it possible to determine, in a rather simple way, a Fourier series representation and a factorization of these functions.
\bs
\bs

\bs

Classification A.M.S. 2010 : 33E05, 11G16, 11J89.
\bs
\bs

{\bf Introduction}
\ms

Since Abel, Jacobi and Weiertsrass, the bases of the theory of the elliptic functions have been well established. Abel proceeds by a length inversion of arc of ellipse (which justifies the denomination of these functions), Jacobi introduces theta functions to define the elliptic functions, finally Weierstrass defines those using series which deliver in an immediate way the double periodicity of those functions.
\ms

Following Weierstrass, the modern presentation of the elliptic functions rest on two pillars: the meromorphic character and the double periodicity. This last property naturally involves well-known representations of the trigonometrical type for Jacobi's elliptic functions - Cf. [Halp], [Whit] et [Armi] - but not for  Weierstrass's elliptic functions.
\ms

We propose in this article a unified framework to obtain, by the same formalism, the elliptic functions of Weierstrass and Jacobi. According to the ideas developed in [Feau], a calculation of contour integral will give a direct access to the Fourier series representations of these elliptic functions. We will also deduce integral representations as well as factorizations of these renormalized elliptic functions.
\bs\ms

{\bf I -- Study of a typical case}
\ms

In order to justify the whole of the study, let us start with a typical case.

As we will establish it with part III, the function $H(s) = \ds \frac{s}{\sin(s)\sin(is)}$ decrease rather quickly ad infinitum when its poles are avoided.

Precisely, for sequences of contours $\gamma_n$, which end up including $\C$, one will establish
\[
\forall z\in ]-\frac{1}{2},\frac{1}{2}[+i]-\frac{1}{2},\frac{1}{2}[, \ds \ \ \lim_{n\to +\infty} \int_{\gamma_n} \frac{s}{\sin(s)\sin(is)}e^{2isz}ds = 0.
\]

The poles of $H$ are localized at $0$, $\pi \Z^*$ and $i \pi \Z^*$, and an easy calculation of residue gives
\[
\forall z\in ]-\frac{1}{2},\frac{1}{2}[+i]-\frac{1}{2},\frac{1}{2}[, \ \ \ds \frac{1}{4\pi} + \sum_{n=1}^{+\infty} (-1)^nn\frac{\cos(2n\pi z)}{\sh(n\pi)} = -\frac{1}{4\pi} - \sum_{n=1}^{+\infty} (-1)^nn\frac{\ch(2n\pi z)}{\sh(n\pi)}.
\]


The $1$-periodic and $i$-periodic character of the two members of this equality makes possible a holomorphic prolongation of the function
\[
\tilde{\wp}(z) = \ds \pi + 4\pi^2\sum_{n=1}^{+\infty} (-1)^nn\frac{\cos(2n\pi z)}{\sh(n\pi)}
\]

on $\abs{Re(z)} < \frac{1}{2}$, $\abs{Im(z)} < \frac{1}{2}$ then on $\C - (\frac{1+i}{2} + \Z + i\Z)$ and by double periodicity : $\forall (n,m)\in\Z^2$, $\tilde{\wp}(z + n + im) = \tilde{\wp}(z)$.
\ms

By integrating $3$ times $\tilde{\wp}$, the singularity into $\frac{1+i}{2}$ becomes artificial: the order of poles are inferior or equal to $3$. It is thus about an elliptic function (Cf. II for a recall of the definition and the elementary properties of these functions).
\ms

Using periodicity, one deduces the existence of a single pole on the elementary domain $[0,1]+i[0,1]$. The parity of $\tilde{\wp}$ makes it possible to conclude to double poles located at $\Lambda = \frac{1+i}{2} + \Z + i\Z$. It is Weierstrass type elliptic function. The parity and nullity at $0$ of $\tilde{\wp}$ then make it possible to assert the existence of $\alpha\in\C$ (it will be shown that $\alpha = 1$) such as
$$\tilde{\wp}(z) = \alpha\sum_{\omega\in \Lambda} \frac{1}{(z-\omega)^2} - \frac{1}{\omega^2}.$$

A first natural result is thus obtained: the Fourier series representations of a Weierstrass type function, which is not current in the literature. The reason comes certainly from the usual choice to place the pole at $0$. We will discuss thereafter the utility of the translation of $\frac{1+i}{2}$ operated here, it is an essential point.
\bs

{\bf II -- Notations and generalities about the elliptic functions}
\ms

Being given $\omega_1$ and $\omega_2$, two non-zero complex numbers such as $\omega_1/\omega_2 \not \in \R$, and let us denote $\Lambda_{\omega_1,\omega_2} = \omega_1\Z + \omega_2\Z$. The lattices $\Lambda_{\omega_1,\omega_2}$ and $\Lambda_{\omega_1,-\omega_2}$ are identical, and we can impose $Im(\omega_2/\omega_1) >0$, \textsl{which we will do systematically thereafter}.
\ms

We pose $\Delta_{\omega_1,\omega_2} = \omega_1[-\frac{1}{2},\frac{1}{2}[ + \omega_2[-\frac{1}{2},\frac{1}{2}[$, this is the basic elementary domain and any translation also will be named elementary domain. It is noticed that $(a + \Delta_{\omega_1,\omega_2})_{a\in \Lambda_{\omega_1,\omega_2}}$ is a partition of $\C$.
\ms

Finally, we set $\omega_0 = \ds \frac{\omega_1+\omega_2}{2}$ and $\tilde{\Delta}_{\omega_1,\omega_2} = \Delta_{\omega_1,\omega_2}-\{-\omega_0\}$.
\bs

Being given $\omega_1$ and $\omega_2$, two non-zero complex numbers such that $\omega_1/\omega_2 \not \in \R$ (and accordingly $Im(\omega_2/\omega_1) >0$), we say as $f : \C \to \C\cup\{\infty\}$ is an elliptic function based on the lattice $\Lambda_{\omega_1,\omega_2}$ if

\be
\item[(i)]
$f$ is a meromorphic function on $\C$,
\item[(ii)]
$f$ is périodic in two directions, says $\omega_1$ and $\omega_2$.
\ee
\ms

The well-known theory of the elliptic functions gives the following results:

\be
\item[(i)]
If $f$ is elliptic and holomorphic on $\C$, then $f$ is a constant function.

\item[(ii)]
The number of poles and zeros in an elementary domain is finite.
Moreover, the sum of the residues being in an elementary domain with no pole on its border is null: if $(a_i)_{i\in I}$ are the poles of $f$ at $a+\Delta_{\omega_1,\omega_2}$, beyond the border, then:
$$\sum_{i\in I} Res(f,a_i) = 0.$$

\item[(iii)]
The sum of the multiplicities of the poles (counted negatively) and the zeros (counted positively) of $f$ in an elementary domain is null.
\ee

To be brief, $(i)$ comes from Liouville's theorem, $(ii)$ and $(iii)$ are obtained by calculating the integrals of $f$ then $f'/f$ on the border of an elementary domain without pole or zero on its border, those being null by symmetry. One will be able, for example, to consult [Whit] or [Armi] on these points.
\bs

Finally, let us notice that $z\mapsto f(z)$ is an elliptic function on the lattice $\Lambda_{\omega_1,\omega_2}$ if and only if $\ds z\mapsto f\pa{\omega_1 z}$ is an elliptic function on $\Lambda_{1,\tau}$, with $\tau = \ds \frac{\omega_2}{\omega_1}$.
\ms
\bs

{\bf III - Weierstrass's modified elliptic function on  $\Lambda_{\omega_1,\omega_2}$}
\ms

To generalize the results of part I, the function $\ds H(s) = \frac{s}{\sin(\omega_1 s)\sin(\omega_2 s)}$ is considered.
\ms

For $n\in\N$, we define a contour $C_n$ connecting by segments,
and in the trigonometrical direction, the points
$(n+\frac{1}{2})\pi(\omega_1^{-1}-\omega_2^{-1})$, $(n+\frac{1}{2})\pi(-\omega_1^{-1}-\omega_2^{-1})$, $(n+\frac{1}{2})\pi(-\omega_1^{-1}+\omega_2^{-1})$ et $(n+\frac{1}{2})\pi(\omega_1^{-1}+\omega_2^{-1})$.

By this choice, the contour passes exactly ``between'' the poles of $H$ on the axes $\omega_1\R$ and $\omega_2\R$. The four segments are noted respectively $C_{nj}$, $1\ie j \ie 4$ (on figure $\omega_1 = 1$ and $\omega_2 = i$):
\bs
\bs
\bs

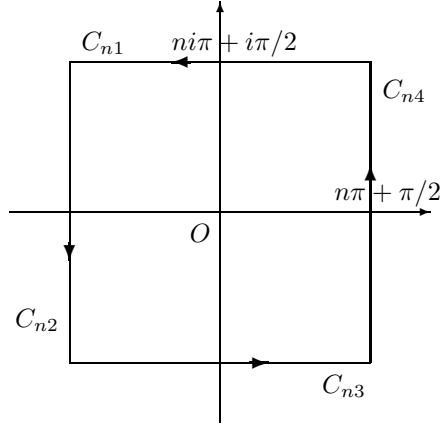
\[
\setlength{\unitlength}{0.8cm}
\begin{picture}(6,4)(-3,-2)

\put(1.9,0.2){$n\pi+\pi/2$} 
\put(-0.8,2.7){$ni\pi+i\pi/2$}

\put(-0.5,-0.5){$O$}

\put(2.7,1.9){$C_{n4}$} 
\put(-2.3,2.7){$C_{n1}$}
\put(1.7,-3){$C_{n3}$}
\put(-3.4,-1.9){$C_{n2}$}

\put(-3.5,0){\vector(1,0){7}} 
\put(0,-3.5){\vector(0,1){7}}

\put(-2.5,-2.5){\line(1,0){5}} 
\put(-2.5,2.5){\line(1,0){5}}
\put(-2.5,-2.5){\line(0,1){5}} 
\put(2.5,-2.5){\line(0,1){5}}

\thicklines
\put(-2.5,-0.5){\vector(0,-1){0.3}} 
\put(2.5,0.5){\vector(0,1){0.3}} 

\put(0.5,-2.5){\vector(1,0){0.3}} 
\put(-0.5,2.5){\vector(-1,0){0.3}} 

\end{picture}
\]
\bs
\bs
\bs

{\bf Theoreme 1} : for $z \in \ ]-\frac{1}{2},\frac{1}{2}[\omega_1 + ]-\frac{1}{2},\frac{1}{2}[\omega_2$, we have $\ds \lim_{n\to +\infty}\int_{C_n} H(s)e^{2isz}ds = 0$.
\bs

\textit{Proof :} Let us verify the result for the segment $C_{n,1}$, the three other cases resulting some by symmetry.

We have to show that $\ds \lim_{n\to +\infty}\int_{[(n+\frac{1}{2})\pi(-1-\tau^{-1}),(n+\frac{1}{2})\pi(1-\tau^{-1})]} \frac{s}{\sin(s)\sin(\tau s)}e^{2isz}ds = 0$, which will give the result by change of variables $\tau = \omega_2/\omega_1$, for $z \in \ ]-\frac{1}{2},\frac{1}{2}[ + ]-\frac{1}{2},\frac{1}{2}[\tau$.
\ms

Let $z  = a + b\tau$, $(a,b)\in \ ]-\frac{1}{2},\frac{1}{2}[^2$.
For $s = u + (n+\frac{1}{2})\pi\tau^{-1}$, with $-(n+\frac{1}{2})\pi \ie u \ie (n+\frac{1}{2})\pi$,

$$\abs{\frac{e^{2isz}}{\sin(s)\sin(\tau s)}} = \abs{\frac{e^{2ias}}{\sin(s)}}\abs{\frac{e^{2ibs\tau}}{\sin(\tau s)}} = \abs{\frac{e^{2ia(n+\frac{1}{2})\pi\tau^{-1}}}{\sin(u+(n+\frac{1}{2})\pi\tau^{-1})}}\abs{\frac{e^{2ibu\tau}}{\cos(\tau u)}}.$$

We denote $\tau^{-1} = \alpha + i\beta$ (let us note that $\beta < 0$) and
\begin{align*}
\ds \abs{\frac{e^{2ia(n+\frac{1}{2})\pi\tau^{-1}}}{\sin(u+(n+\frac{1}{2})\pi\tau^{-1})}} &= \ds \abs{\frac{e^{-2a(n+\frac{1}{2})\pi\beta}}{\sin((u+(n+\frac{1}{2})\pi\alpha) + i(n+\frac{1}{2})\pi\beta)}}\hfill \\
&= \ds  \abs{\frac{e^{-(n+\frac{1}{2})\pi\beta + i(u+(n+\frac{1}{2})\pi\alpha)}}{{\rm \sh}((n+\frac{1}{2})\pi\beta - i(u+(n+\frac{1}{2})\pi\alpha))}}e^{(1-2a)(n+\frac{1}{2})\pi\beta}.\hfill \\
\end{align*}

The function $\ds \mapsto \abs{\frac{e^{v}}{{\rm \sh}(v)}}$ is bounded on $\abs{Re(v)} > 1$, whereas $\ds u\mapsto \abs{\frac{e^{2ibu\tau}}{\cos(\tau u)}}$ is bounded on $\R$, because $\abs{2b} < 1$. Now, we can find $M > 0$ depending only on $z$ and $\tau$ such as
$\ds \abs{\frac{e^{2isz}}{\sin(s)\sin(\tau s)}} \ie M e^{(1-2a)(n+\frac{1}{2})\pi\beta}$, and the result is thus deduced.
\bs\bs

We have to integrate on segments of direction $\omega_1$ and $\omega_2$ passing exactly between the poles of $H$ and to use the method of residues.
\ms

All the pole of $H$ are simple :

\be
\item[(i)]
$s = 0$, and the residue is worth $\ds \frac{1}{\omega_1\omega_2}$.
\item[(ii)]
$s = n\pi/\omega_1$, $n\in\Z^*$, and the residue is worth $\ds\frac{(-1)^n n\pi}{\omega_1^2\sin(n\pi\omega_2/\omega_1)}e^{2in\pi z/\omega_1}$.
\item[(iii)]
$s = n\pi/\omega_2$, $n\in\Z^*$, and the residue is worth $\ds\frac{(-1)^n n\pi}{\omega_2^2\sin(n\pi\omega_1/\omega_2)}e^{2in\pi z/\omega_2}$.
\ee
\ms

The calculus of residues proposed to theorem 1 gives finally: $\forall z \in \ ]-\frac{1}{2},\frac{1}{2}[\omega_1 + ]-\frac{1}{2},\frac{1}{2}[\omega_2$,

$$\frac{\omega_2}{\omega_1}\sum_{n=1}^{+\infty} \frac{(-1)^n n}{\sin(n\pi\omega_2/\omega_1)}\cos(2nz\pi/\omega_1) = - \frac{1}{2\pi} - \frac{\omega_1}{\omega_2}\sum_{n=1}^{+\infty} \frac{(-1)^n n}{\sin(n\pi\omega_1/\omega_2)}\cos(2nz\pi/\omega_2).$$
\ms

Let $\Theta(z)$ be this common value, the periodic character of each member of the equality makes possible a holomorphic prolongation of $\Theta$ on $\abs{Im(z/\omega_1)} < \frac{1}{2} Im(\omega_2/\omega_1)$, $\abs{Im(z/\omega_2)} < -\frac{1}{2} Im(\omega_1/\omega_2)$ then $\C - (\omega_0 + \omega_1\Z + \omega_2\Z)$ by double-periodicity : $\Theta(z + \omega_1 n + \omega_2 m) = \Theta(z)$, $(n,m)\in\Z^2$.

After three integrations, it is noted that the poles of $\Theta$ are all doubles with null residues and placed on the relocated lattice $\Lambda^{\star}_{\omega_1,\omega_2} = \omega_0 + \Lambda_{\omega_1,\omega_2}$, with $\omega_0 = \frac{\omega_1+\omega_2}{2}$.
\ms

According to Liouville's theorem, except for a multiplicative constant and an additive constant, it is about

\begin{align*}
\tilde{\wp}_{\omega_1,\omega_2}(z) &= \ds \sum_{a\in \Lambda^{\star}_{\omega_1,\omega_2}} \frac{1}{(z-a)^2} - \frac{1}{a^2}\hfill \\
 &= \ds \sum_{(n,m)\in\Z} \frac{1}{(z+(n+\frac{1}{2})\omega_1+(m+\frac{1}{2})\omega_2)^2} - \frac{1}{((n+\frac{1}{2})\omega_1+(m+\frac{1}{2})\omega_2)^2}.\hfill \\
\end{align*}

The meromorphic function $\tilde{\wp}_{\omega_1,\omega_2}$ presents double poles with null residues on $\Lambda^{\star}_{\omega_1,\omega_2}$, as well as double roots on $\Lambda_{\omega_1,\omega_2}$. This is the elliptic function which appears naturally here and which is studied now.
\bs

\textsc{Fourier series representations of $\tilde{\wp}_{\omega_1,\omega_2}$}
\ms

There thus exists $\alpha\in \C^*$ such as $\tilde{\wp}_{\omega_1,\omega_2}(z) = \alpha(\Theta(z) - \Theta(0))$.

We will find $\alpha$ while studying  $\ds \Theta(z) = \frac{1}{\alpha(z-\omega_0)^2} + O(1)$ in the vicinity of $\omega_0$.

While developing, it is noted that $\ds \abs{(-1)^n\frac{\cos(2n\pi z/\omega_1)}{\sin(n\pi\omega_2/\omega_1)} + i e^{-2in\pi (z-\omega_0)/\omega_1}}e^{n\pi Im(\omega_2/\omega_1)}$ is bounded independently of $n\in\N^*$ and $z\in [0,\frac{1}{2}[\omega_1 + [0,\frac{1}{2}[\omega_2$. If we introduce $\psi(z) = \ds  -i\frac{\omega_2}{\omega_1}\sum_{n=1}^{+\infty} ne^{-2in\pi(z-\omega_0)/\omega_1}$, the function $\Theta - \psi$ is bounded on  $[0,\frac{1}{2}[\omega_1 + [0,\frac{1}{2}]\omega_2$.
\ms

A direct calculation shows that $\psi(z) = \ds -i\frac{\omega_2}{\omega_1} \frac{e^{-2i\pi(z -\omega_0)/\omega_1}}{(1-e^{-2i\pi(z -\omega_0)/\omega_1})^2}$ and $\ds \lim_{z\to \omega_0} \pa{z - \omega_0}^2\psi(z) = \frac{i\omega_1\omega_2}{4\pi^2}$.

We deduce $\ds \tilde{\wp}(z) = -i\frac{4\pi^2}{\omega_1\omega_2}(\Theta(z) - \Theta(0))$.
\ms

Thus, the function $\tilde{\wp}_{\omega_1,\omega_2}$ is even, admits a double pole  at $\omega_0$, cancels itself twice at $0$ and it comes:

\begin{align}
\ds \tilde{\wp}_{\omega_1,\omega_2}(z)&= \ds -\frac{i4\pi^2}{\omega_1^2}\pa{\sum_{n=1}^{+\infty} \frac{(-1)^n n}{\sin(n\pi\omega_2/\omega_1)}\cos(2n\pi \frac{z}{\omega_1}) - \sum_{n=1}^{+\infty} \frac{(-1)^n n}{\sin(n\pi\omega_2/\omega_1)}}\hfill \\
 &= \ds \frac{i4\pi^2}{\omega_2^2}\pa{\sum_{n=1}^{+\infty} \frac{(-1)^n n}{\sin(n\pi\omega_1/\omega_2)}\cos(2n\pi \frac{z}{\omega_2}) - \sum_{n=1}^{+\infty} \frac{(-1)^n n}{\sin(n\pi\omega_1/\omega_2)}}\hfill \\
  &= \ds \frac{\pi^2}{\omega_1^2}\sum_{n\in\Z} \pa{\frac{1}{\sin((z -\omega_0-\omega_2 n)\pi/\omega_1)^2} - \frac{1}{\sin((\omega_0+\omega_2 n)\pi/\omega_1)^2}}.\hfill \\
  &= \ds \frac{\pi^2}{\omega_2^2}\sum_{n\in\Z} \pa{\frac{1}{\sin((z -\omega_0-\omega_1 n)\pi/\omega_2)^2} - \frac{1}{\sin((\omega_0+\omega_1 n)\pi/\omega_2)^2}}.\hfill
\end{align}

The first two equalities are valid for $\abs{Im(z/\omega_1)} < \frac{1}{2} Im(\omega_2/\omega_1)$ and $\abs{Im(z/\omega_2)} < -\frac{1}{2}Im(\omega_1/\omega_2)$ respectively. The two last ones come from Euler's equality $\sin(z)^{-2} = \ds \sum_{n\in\Z} \frac{1}{(z-n\pi)^2}$ and are valid on $\C-\Lambda^{\star}_{\omega_1,\omega_2}$.
\bs

It should be noted that in the typical case studied with part I, the symmetry of the function gives an explicit formula for the constant: $\ds \sum_{n=1}^{+\infty} \frac{(-1)^n n}{\sin(in\pi)} = -i\sum_{n=1}^{+\infty} \frac{(-1)^n n}{\sh(n\pi)} = \frac{-i}{4\pi}$.
\bs

\textsc{Link with Weierstrass's elliptic function}
\ms

One points out the definition of Weiertrass's function :
\begin{equation}
\wp_{\omega_1,\omega_2}(z) = \frac{1}{z^2} + \sum_{(n,m)\not=(0,0)} \frac{1}{(z+n\omega_1+m\omega_2)^2} - \frac{1}{(n\omega_1+m\omega_2)^2}.\label{Wp}
\end{equation}

It comes immediately:
\[\tilde{\wp}_{\omega_1,\omega_2}(z) = \wp_{\omega_1,\omega_2}(z+ \omega_0) - \wp_{\omega_1,\omega_2}(\omega_0).\]

Let us note the differences:

\be
\item[(i)]
Poles of $\wp_{\omega_1,\omega_2}$ are localised on $\Lambda_{\omega_1,\omega_2}$ whereas those of $\tilde{\wp}_{\omega_1,\omega_2}$ are on $\Lambda_{\omega_1,\omega_2}^\star$.
\item[(ii)]
Zeros of $\tilde{\wp}_{\omega_1,\omega_2}$ are localised on $\Lambda_{\omega_1,\omega_2}$, whereas those of $\wp_{\omega_1,\omega_2}$ are unknown.
\item[(iii)]To define $\wp_{\omega_1,\omega_2}$, one particularizes the pole at $0$ which disymmetrizes the function whereas $\tilde{\wp}_{\omega_1,\omega_2}$ is symmetrical by construction.
\item[(iv)]
The choice of normalization of $\wp_{\omega_1,\omega_2}$ is closely related to the parameterization of a specifical family of the cubic ones. We will see that the choice retained for $\tilde{\wp}_{\omega_1,\omega_2}$ allows to parameterize another family of cubic curves, important also.
\ee

\ms
\bs

{\bf IV - Weiertrass functions on the reduced lattice $\Lambda_{\tau} = \Z + \tau\Z$}
\ms

As previously noted, the study of the elliptic functions on the lattice $\Lambda_{\omega_1,\omega_2}$ is brought back under investigation on the lattice $\Lambda_{1,\tau} = \Lambda_\tau = \Z + \tau\Z$. This will be largely exploited for the study of modular forms.

Let us summarize knowledge concerning the functions $\wp_\tau = \wp_{1,\tau}$ and $\tilde{\wp}_\tau = \tilde{\wp}_{1,\tau}$ on lattice $\Lambda_{\tau}$.

\begin{equation}
\wp_\tau(z) = \frac{1}{z^2} + \sum_{(n,m)\not=(0,0)} \frac{1}{(z+n+m\tau)^2} - \frac{1}{(n+m\tau)^2}, \ \ \forall z\in \C-\Lambda_\tau \label{Wp1}
\end{equation}
and
\begin{equation}
\tilde{\wp}_\tau(z) = \sum_{(n,m)\in\Z^2} \frac{1}{(z+(n+\frac{1}{2})+(m+\frac{1}{2})\tau)^2} - \frac{1}{((n+\frac{1}{2})+(m+\frac{1}{2})\tau)^2}, \ \ \forall z\in \C-(\frac{1+\tau}{2}+\Lambda_\tau). \label{Wpt}
\end{equation}

These functions are connected by $\tilde{\wp}_\tau(z) = \wp_\tau(z+ \frac{1}{2} + \frac{\tau}{2}) - \wp_\tau(\frac{1}{2}+\frac{\tau}{2})$.
\bs

\textsc{$\bullet$ The $\tilde{\wp}$ function}
\be
\item[-] The integral of the functional equation:
$\forall z \in \ ]-\frac{1}{2},\frac{1}{2}[ + ]-\frac{1}{2},\frac{1}{2}[\tau$,

$$\int_{C_\infty} \frac{s}{\sin(s)\sin(\tau s)}e^{2isz}ds = 0.$$

\item[-] The calculus of residues:
$\forall z \in \ ]-\frac{1}{2},\frac{1}{2}[ + ]-\frac{1}{2},\frac{1}{2}[\tau$,

\[\frac{1}{2\pi} + \tau\sum_{n=1}^{+\infty} \frac{(-1)^n n}{\sin(n\pi\tau)}\cos(2nz\pi) + \frac{1}{\tau}\sum_{n=1}^{+\infty} \frac{(-1)^n n}{\sin(n\pi/\tau)}\cos(2nz\pi/\tau) = 0.
\]

\item[-] The Fourier series representations: $\forall z = \alpha + \beta\tau$, $(\alpha,\beta)\in\R^2$,

\begin{align}
\ds \tilde{\wp}_\tau(z) &= \ds 4i\pi^2\pa{\sum_{n=1}^{+\infty} \frac{(-1)^n n}{\sin(n\pi\tau)} - \sum_{n=1}^{+\infty} \frac{(-1)^n n}{\sin(n\pi\tau)}\cos(2n\pi z)}\hfill & \ {\rm when} \ \abs{\beta} < \frac{1}{2}\\
  &= \ds \frac{4i\pi^2}{\tau^2}\pa{-\sum_{n=1}^{+\infty} \frac{(-1)^n n}{\sin(n\pi/\tau)} + \sum_{n=1}^{+\infty} \frac{(-1)^n n}{\sin(n\pi/\tau)}\cos(2n\pi z/\tau)}\hfill & \ {\rm when} \ \abs{\alpha} < \frac{1}{2}.
\end{align}

\item[-] The integral representations: $\forall z = \alpha + \beta\tau$, $(\alpha,\beta)\in\R^2$,

\begin{align*}
\tilde{\wp}_{\tau}(z) &= \ds -8i\pi^2 \int_{C_{1}} \frac{s}{\sin(s)\sin(\tau s)} (\cos(2sz)-1)ds \hfill & \ {\rm when} \ \abs{\beta} < \frac{1}{2}\\
 &= \ds 8i\pi^2 \int_{C_{\tau}} \frac{s}{\sin(s)\sin(\tau s)} (\cos(2sz)-1)ds\hfill & \ \ {\rm when} \ \abs{\alpha} < \frac{1}{2}
\end{align*}

where $C_1$ is a directed contour which includes the poles $n\pi$, $n\in\N^*$, and $C_\tau$ a directed contour which includes the poles $n\frac{\pi}{\tau}$, $n\in\N^*$.
\ee
\bs

\textsc{$\bullet$ La fonction $\wp$}
\be
\item[-] The integral of the functional equation:
$\forall z \in \ ]0,1[ + ]0,1[\tau$,

$$\int_{C_\infty} \frac{s}{\sin(s)\sin(\tau s)}e^{is(2z-1-\tau)}ds = 0.$$

\item[-] The calculus of residues:
$\forall z \in \ ]0,1[ + ]0,1[\tau$,

$$\frac{1}{2\pi} + \tau\sum_{n=1}^{+\infty} \frac{n}{\sin(n\pi\tau)}\cos(n\pi(2z-\tau)) + \frac{1}{\tau}\sum_{n=1}^{+\infty} \frac{n}{\sin(n\pi/\tau)}\cos(n(2z-1)\pi/\tau) = 0.$$

\item[-] The Fourier series representations:
$\forall z = \alpha + \beta\tau$, $(\alpha,\beta)\in\R^2$,

\begin{align*}
\ds \wp_\tau(z) &= \ds c_\tau - 4i\pi^2 \sum_{n=1}^{+\infty} \frac{n}{\sin(n\pi\tau)}\cos(n\pi (2z-\tau))\hfill & \ {\rm when} \ \beta\in ]0,1[ \\
  &= \ds d_\tau + \frac{4i\pi^2}{\tau^2}\sum_{n=1}^{+\infty} \frac{n}{\sin(n\pi/\tau)}\cos(n(2z-1)\pi/\tau)\hfill & \ {\rm when} \ \alpha\in ]0,1[.
\end{align*}

The terms $c_\tau$ and $d_\tau$ depending only on $\tau$ and satisfying $c_\tau - d_\tau = -\frac{2}{\tau}$.

\item[-] The integral representations:
$\forall z = \alpha + \beta\tau$, $(\alpha,\beta)\in\R^2$,

\begin{align*}
\wp_{\tau}(z) &= \ds c_\tau - 8i\pi^2 \int_{C_{1}} \frac{s}{\sin(s)\sin(\tau s)} \cos(s(2z-1-\tau))ds \hfill & \ {\rm when} \ \beta\in ]0,1[ \\
 &= \ds d_\tau + 8i\pi^2 \int_{C_{1}} \frac{s}{\sin(s)\sin(\tau s)} \cos(s(2z-1-\tau))ds.\hfill & \ {\rm when} \ \alpha\in ]0,1[
\end{align*}

where $C_1$ is a directed contour which includes the poles $n\pi$, $n\in\N^*$, and $C_\tau$ a directed contour which includes the poles $n\frac{\pi}{\tau}$, $n\in\N^*$.
\ee
\bs
\ms

{\bf V - Jacobi's elliptic functions}
\ms

Historically preceding Weieirstrass's theory, the Jacobi elliptic functions constitute the other standard family of elliptic functions.

We return to [Armi] and [Whit] for the definitions of these functions. Let us note simply that the periods are related to two real numbers $K$ and $K'$ (defined themselves by a parameter $k$), that one pose $\tau = iK'/K$ and $q = e^{i\pi\tau}$ (caution: it is not the usual convention in modular forms theory). To establish the link with what follows, one finds the Fourier series representations of $ {\rm Sn} $, $ {\rm Cn} $ and $ {\rm dn} $ (with the notation in $q$) in [Armi] and [Whit].

The approach developed in part III with contour calculus $\ds \lim_{n\to+\infty} \int_{C_n} H(s)e^{2isz}ds = 0$ remains valid to obtain the functions of Jacobi. The demonstrations of convergence are similar to those developed partly III and will be omitted. In all the cases, after a transformation $z \mapsto \alpha z$ where $\alpha \in \C$, the double periodicity is reduced to the periods $ (1,2\tau)$ or $ (2,2\tau)$ with $Im(\tau) > 0$. Following the example of the choice of lattice $\Lambda_\tau$ in the previous section, the reduced form will be very useful for the study of modular forms.
\bs

\textsc{$\bullet$ Jacobi's ${\rm dn}$ function}
\be
\item[-] The integral of the functional equation:
$\forall z \in \ \frac{1}{2} + ]-\frac{1}{2},\frac{1}{2}[ + ]-\frac{1}{2},\frac{1}{2}[\tau$,

$$\int_{C_\infty} \frac{1}{\sin(s)\cos(\tau s)}e^{is(2z-1)}ds = 0.$$

\item[-] The calculus of residues:
$\forall z \in \ \frac{1}{2} + ]-\frac{1}{2},\frac{1}{2}[ + ]-\frac{1}{2},\frac{1}{2}[\tau$,

$$1 + 2\sum_{n=1}^{+\infty} \frac{\cos(2n\pi z)}{\cos(n\pi\tau)} - \frac{2}{\tau}\sum_{n=0}^{+\infty} \frac{(-1)^n \cos((n+\frac{1}{2})\pi(2z-1)/\tau)}{\sin((n+\frac{1}{2})\pi/\tau)} = 0.$$

\item[-] The Fourier series representations:
$\forall z = \frac{1}{2} + \alpha + \beta\tau$, $(\alpha,\beta)\in\R^2$,

\begin{align}
\ds DN(z,\tau) = \frac{2K}{\pi} \ {\rm dn}(2Kz) &= \ds 1 + 4\sum_{n=1}^{+\infty} \frac{q^{n}\cos(2n\pi z)}{1 + q^{2n}}\hfill & \ {\rm when} \ \abs{\beta} < \frac{1}{2}  \\
 &= 1 + 2\sum_{n=1}^{+\infty} \frac{\cos(2n\pi z)}{\cos(n\pi\tau)}\hfill & \ {\rm when} \ \abs{\beta} < \frac{1}{2} \\
  &= \ds \frac{2}{\tau}\sum_{n=0}^{+\infty} \frac{(-1)^n \cos((n+\frac{1}{2})\pi(2z-1)/\tau)}{\sin((n+\frac{1}{2})\pi/\tau)}\hfill & \ {\rm when} \ \abs{\alpha} < \frac{1}{2}.
\end{align}

\item[-] The integral representations:
$\forall z = \frac{1}{2} + \alpha + \beta\tau$, $(\alpha,\beta)\in\R^2$,

\begin{align*}
\ds DN(z,\tau) = \frac{2K}{\pi} \ {\rm dn}(2Kz) &= \ds 1 + 4i\pi \int_{C_{1}} \frac{1}{\sin(s)\cos(\tau s)}\cos(s(2z-1))ds \hfill & \ {\rm when} \ \abs{\beta} < \frac{1}{2} \\
 &= \ds -4i\pi\int_{C_{\tau}} \frac{1}{\sin(s)\cos(\tau s)}\cos(s(2z-1))ds\hfill & \ {\rm when} \ \abs{\alpha} < \frac{1}{2},
\end{align*}

where $C_1$ includes the poles $n\pi$, $n\in\N$, and $C_\tau$ the poles $(n+\frac{1}{2})\frac{\pi}{\tau}$.
\ee
\bs

\textsc{$\bullet$ Jacobi's ${\rm sn}$ function}
\be
\item[-] The integral of the functional equation:
$\forall z \in \ \frac{1}{2} + ]-\frac{1}{2},\frac{1}{2}[ + ]-\frac{1}{2},\frac{1}{2}[\tau$,

$$\int_{C_\infty} \frac{1}{\cos(s)\sin(\tau s)}e^{is(2z-1)}ds = 0.$$

\item[-] The calculus of residues:
$\forall z \in \ \frac{1}{2} + ]-\frac{1}{2},\frac{1}{2}[ + ]-\frac{1}{2},\frac{1}{2}[\tau$,

$$1 + 2\sum_{n=1}^{+\infty} \frac{(-1)^n \cos(n\pi(2z-1)/\tau)}{\cos(n\pi/\tau)} - 2\tau\sum_{n=0}^{+\infty} \frac{\sin((2n+1)\pi z)}{\sin((n+\frac{1}{2})\pi\tau)} = 0.$$

\item[-] The Fourier series representations:
$\forall z = \frac{1}{2} + \alpha + \beta\tau$, $(\alpha,\beta)\in\R^2$,

\begin{align}
\ds SN(z,\tau) = \frac{Kk}{2\pi} \ {\rm sn}(2Kz) &= \ds \sum_{n=0}^{+\infty} \frac{q^{n+\frac{1}{2}}\sin((2n+1)\pi z)}{1 - q^{2n+1}}\hfill & \ {\rm when} \ \abs{\beta} < \frac{1}{2} \\
 &= \frac{i}{2}\sum_{n=0}^{+\infty} \frac{\sin((2n+1)\pi z)}{\sin((n+\frac{1}{2})\pi\tau)}\hfill & \ {\rm when} \ \abs{\beta} < \frac{1}{2} \\
  &= \ds \frac{i}{4\tau} + \frac{i}{2\tau}\sum_{n=1}^{+\infty} \frac{(-1)^n \cos(n\pi(2z-1)/\tau)}{\cos(n\pi/\tau)}\hfill & \ {\rm when} \ \abs{\alpha} < \frac{1}{2}.
\end{align}

\item[-] The integral representations:
$\forall z = \frac{1}{2} + \alpha + \beta\tau$, $(\alpha,\beta)\in\R^2$,

\begin{align*}
\ds SN(z,\tau) = \frac{Kk}{2\pi} \ {\rm sn}(2Kz) &= \ds \pi \int_{C_{1}} \frac{1}{\cos(s)\sin(\tau s)}\cos(s(2z-1))ds \hfill & \ {\rm when} \ \abs{\beta} < \frac{1}{2} \\
 &= \ds \frac{i}{4\tau} - \pi\int_{C_{\tau}} \frac{1}{\cos(s)\sin(\tau s)}\cos(s(2z-1))ds\hfill & \ {\rm when} \ \abs{\alpha} < \frac{1}{2},
 \end{align*}

where $C_1$ includes the poles $(n+\frac{1}{2})\pi$, $n\in\N$, and $C_\tau$ the poles $n\pi/\tau$, $n\in\N^*$.
\ee
\bs

\textsc{$\bullet$ Jacobi's ${\rm cn}$ function}
\be
\item[-] The integral of the functional equation:
$\forall z \in \ \frac{1}{2} + ]-\frac{1}{2},\frac{1}{2}[ + ]-\frac{1}{2},\frac{1}{2}[\tau$,

$$\int_{C_\infty} \frac{1}{\cos(s)\cos(\tau s)}e^{is(2z-1)}ds = 0.$$

\item[-] The calculus of residues:
$\forall z \in \ \frac{1}{2} + ]-\frac{1}{2},\frac{1}{2}[ + ]-\frac{1}{2},\frac{1}{2}[\tau$,

$$\sum_{n=1}^{+\infty} \frac{(-1)^n \sin((n+\frac{1}{2})\pi(2z-1)/\tau)}{\cos((n+\frac{1}{2})\pi/\tau)} - \tau\sum_{n=0}^{+\infty} \frac{\cos((2n+1)\pi z)}{\cos((n+\frac{1}{2})\pi\tau)} = 0.$$

\item[-] The Fourier series representations:
$\forall z = \frac{1}{2} + \alpha + \beta\tau$, $(\alpha,\beta)\in\R^2$,

\begin{align}
\ds CN(z,\tau) = \frac{Kk}{2\pi} \ {\rm cn}(2Kz) &= \ds \sum_{n=0}^{+\infty} \frac{q^{n+\frac{1}{2}}\cos((2n+1)\pi z)}{1 + q^{2n+1}}\\
 &= \ds \frac{1}{2}\sum_{n=0}^{+\infty} \frac{\cos((2n+1)\pi z)}{\cos((n+\frac{1}{2})\pi\tau)}\hfill & \ {\rm when} \ \abs{\beta} < \frac{1}{2} \\
  &= \ds \frac{1}{2\tau}\sum_{n=0}^{+\infty} \frac{(-1)^n \sin((n+\frac{1}{2})\pi(2z-1)/\tau)}{\cos((n+\frac{1}{2})\pi/\tau)}\hfill & \ {\rm when} \ \abs{\alpha} < \frac{1}{2}.
\end{align}

\item[-] The integral representations:
$\forall z = \frac{1}{2} + \alpha + \beta\tau$, $(\alpha,\beta)\in\R^2$,

\begin{align*}
\ds CN(z,\tau) = \frac{Kk}{2\pi} \ {\rm cn}(2Kz) &= \ds i\pi \int_{C_{1}} \frac{1}{\cos(s)\cos(\tau s)}\sin(s(2z-1))ds \hfill & \ {\rm when} \ \abs{\beta} < \frac{1}{2} \\
 &= \ds -i\pi\int_{C_{\tau}} \frac{1}{\cos(s)\cos(\tau s)}\sin(s(2z-1))ds\hfill & \ {\rm when} \ \abs{\alpha} < \frac{1}{2}.
\end{align*}

where $C_1$ includes the poles $(n+\frac{1}{2})\pi$, $n\in\N$, and $C_\tau$ the poles $(n+\frac{1}{2})\pi/\tau$, $n\in\N^*$.

\ms

\ee
\bs

\textsc{$\bullet$ Jacobi's ${\rm nd}$, ${\rm cd}$ and ${\rm sd}$ functions}
\ms

The reader will verify without difficulty that one obtains in the same way Fourier series representations and integral representations for the other Jacobi's elliptic functions. For an application to the modular forms, the renormalized functions ${\rm nd}$, ${\rm cd}$ and ${\rm sd}$ will be useful. One considers respectively:

$$\int_{C_\infty} \frac{1}{\sin(s)\cos(\tau s)}e^{2isz}ds = 0, \ \int_{C_\infty} \frac{1}{\cos(s)\sin(\tau s)}e^{2isz}ds = 0 \ {\rm et} \ \int_{C_\infty} \frac{1}{\cos(s)\cos(\tau s)}e^{2isz}ds = 0,$$

i.e., by changing $2z-1$ by $2z$ in the previous integrals. As previously, the following Fourier series representations are valid on $z = \alpha + \beta \tau$ for $\abs{\beta} < \frac{1}{2}$ then for $\abs{\alpha} < \frac{1}{2}$. It comes:

\begin{align}
\ds ND(z,\tau) = \ds \frac{2k'K}{\pi} \ {\rm nd}(2Kz) &= \ds 1 + 4\sum_{n=1}^{+\infty} \frac{(-1)^nq^{n}\cos(2n\pi z)}{1 + q^{2n}} = 1 + 2\sum_{n=1}^{+\infty} \frac{(-1)^n\cos(2n\pi z)}{\cos(n\pi\tau)}\hfill \\
 &= \ds \frac{2}{\tau}\sum_{n=0}^{+\infty} \frac{(-1)^n \cos((2n+1)\pi z/\tau)}{\sin((n+\frac{1}{2})\pi/\tau)}\hfill \\
\ds CD(z,\tau) = \ds \frac{Kk}{2\pi} \ {\rm cd}(2Kz) &= \ds \sum_{n=0}^{+\infty} \frac{(-1)^nq^{n+\frac{1}{2}}\cos((2n+1)\pi z)}{1 - q^{2n+1}} = \frac{i}{2}\sum_{n=0}^{+\infty} \frac{(-1)^n\cos((2n+1)\pi z)}{\sin((n+\frac{1}{2})\pi\tau)}\hfill\hfill \\
 &= \ds \frac{i}{4\tau} + \frac{i}{2\tau}\sum_{n=1}^{+\infty} \frac{(-1)^n \cos(2n\pi z/\tau)}{\cos(n\pi/\tau)}\hfill \\
\ds SD(z,\tau) = \ds \frac{Kkk'}{2\pi} \ {\rm sd}(2Kz) &= \ds \sum_{n=0}^{+\infty} \frac{(-1)^nq^{n+\frac{1}{2}}\sin((2n+1)\pi z)}{1 + q^{2n+1}} = \ds \frac{1}{2}\sum_{n=0}^{+\infty} \frac{(-1)^n\sin((2n+1)\pi z)}{\cos((n+\frac{1}{2})\pi\tau)}\hfill \\
&= \ds -\frac{1}{2\tau}\sum_{n=0}^{+\infty} \frac{(-1)^n \sin((2n+1)\pi z/\tau)}{\cos((n+\frac{1}{2})\pi/\tau)}\hfill
\end{align}

\ms
\bs

{\bf V - Factorization of elliptic functions}
\ms

The zeros and the poles of the elliptic functions studied previously are quite localised which will make it possible to factorize these functions. The results presented here will have a significant impact in the studies of the modular forms which result from it. We start by factorizing $\twp_\tau$ then the Jacobi's functions, these functions being renormalized from the point of view to come.
\ms

\textsc{$\bullet$ Factorization of $\ds \twpa(z,\tau) = -\frac{1}{16\pi^2}\twp_\tau(z)$}
\ms

In order to simplify the expressions in the modular equalities, one standardizes $\twp_\tau$ by also highlighting the variable $\tau$:
\[\twpa(z,\tau) = -\frac{1}{16\pi^2}\twp_\tau(z).\]

One points out the definition of the reduced lattice $\Lambda_\tau = \Z + \tau \Z$ and one introduces $\Lambda_\tau^\star = \frac{1+\tau}{2} +\Z + \tau\Z$.
\ms

The variable $\tau$ being fixed, in order to compensate zeros and poles of $z\mapsto \twpa(z,\tau)$, the following meromorphic function is introduced:
\[F(z,\tau) = \prod_n\prod_m\pa{\frac{z+m+n\tau}{z+\omega_0+m+n\tau}}^2.\]

While gathering according to the variable $m$ and while symmetrizing in order to have terms converging, the well known Euler's relation $\ds \sin(\pi z) = \pi z\lim_{N\to+\infty} \prod_{n=-N, \ n\not=0}^N \pa{\frac{z+n}{n}}$ leads us to:
\ms

\begin{center}
\begin{tabular}{lcl}
$\ds F(z,\tau)$ & $=$ & $\ds \lim_{N\to +\infty} \prod _{n=-N}^{N}{\frac { \sin \left( \pi \, \left( z+n\tau
 \right)  \right)  ^{2}}{\sin \left( \pi \, \left( z-1/2-1/2\,
\tau+n\tau \right)  \right) \sin \left( \pi \, \left( z+1/2+1/2\,\tau+n\tau
 \right)  \right) }}$\\
 & $=$ & $\ds \prod _{n=-\infty}^{+\infty}{\frac { \sin \left( \pi \, \left( z+n\tau
 \right)  \right)  ^{2}}{\sin \left( \pi \, \left( z-1/2-1/2\,
\tau+n\tau \right)  \right) \sin \left( \pi \, \left( z+1/2+1/2\,\tau+n\tau
 \right)  \right) }}$\\
\end{tabular}
\end{center}
\ms

which, always to $\tau$ fixed, converges uniformly on very compact of $\C - \Lambda_\tau^\star$ according to the variable $z$ with poles of order two. The limit is thus meromorphic, $1$ and $\tau$ periodic.
\ms

Thanks to Liouville's theorem, it is possible to affirm that $z\mapsto \twpa(z,\tau)F(z,\tau)^{-1}$ is a constant in variable $z$ and there exists a function $\tau \mapsto a(\tau) $ such as:
\[\forall (z,\tau)\in (\C-\Lambda^\ast)\times {\cal H}, \  \ \ \twpa(z,\tau) =
a(\tau) \prod _{n=-\infty}^{+\infty}{\frac {\sin \left( \pi \, \left( z+nt
 \right)  \right) ^{2}}{\sin \left( \pi \, \left( z-1/2-1/2\,
t+nt \right)  \right) \sin \left( \pi \, \left( z+1/2+1/2\,t+nt
 \right)  \right) }}.\]

For the continuation, and in accordance with the use for the modular forms, one poses $q = e^{2i\pi\tau}$ and one seeks $\alpha(q) = a(\tau)$ satisfying:

\[\twpa(z,\tau) = -\alpha(q) \frac{q^{1/2}}{4 \sin \left( \pi \,z \right) ^{2}} \left( \prod _{n=0}^{+\infty}{\frac { \left( 1-{{\rm e}^{2\,
i\pi \,z}}{q}^{n} \right)  \left( 1-{{\rm e}^{-2\,i\pi \,z}}{q}^{n}
 \right) }{ \left( 1+{{\rm e}^{2\,i\pi \,z}}{q}^{n+1/2} \right) 
 \left( 1+{{\rm e}^{-2\,i\pi \,z}}{q}^{n+1/2} \right) }} \right) ^{2}.\]

From the Fourier series representation:
\[
\twpa(\frac{1}{2},2\tau) = \ds \frac{i}{2}\sum_{n=0}^{+\infty} \frac{2n+1} {\sin(2(2n+1)\pi\tau)} \label{Wpt2}
\]

one obtains easily $\ds \twpa(\frac{1}{2},\tau) = \sum_{n=0}^{+\infty} \sigma_1(2n+1) q^{2n+1}$ where $\sigma_1(n) = \ds \sum_{d\mid n, \ d>0} d$.

We recall a modular equality:
\[\forall \tau\in {\cal H}, \ \ \frac{\eta(4\tau)^8}{\eta(2\tau)^4} = \sum_{n=0}^{+\infty} \sigma_1(2n+1) q^{2n+1},\]
where $\eta(\tau) = \ds q^{-1/24}\prod_{n=1}^{+\infty} (1-q^n)$ is Dedekind's $\eta$ function, see [Rou] page 2 or [Koh] page 146 for instance. And consequently,

\[\twpa(\frac{1}{2},\tau) =
q^{1/2} \prod _{k=0}^{+\infty}{\frac { \left( 1-{q}^{2\,k+2} \right) ^{4}}{
 \left( 1-{q}^{2\,k+1} \right) ^{4}}}.
\]

Thanks to the equality $\ds \prod_{n=0}^{+\infty}(1+q^n). \prod_{n=0}^{+\infty}(1-q^{2n+1}) = 2$, we find

\[\twpa(\frac{1}{2},\tau) = -\alpha(q)\frac{q^{1/2}}{4} \prod _{n=0}^{+\infty}{\frac { \left( 1+{q}^{n} \right) ^{4}}
{ \left( 1-{q}^{n+1/2} \right) ^{4}}} = -4\alpha(q)q^{1/2} \prod _{n=0}^{+\infty}{\frac {1}
{ \left( 1-{q}^{2n+1} \right) ^{4} \left( 1-{q}^{n+1/2} \right) ^{4}}}\]

and $-4\alpha(q) = \ds \prod_{n=0}^{+\infty} (1-q^{2n+2})^4(1-q^{n+1/2})^4$. Finally,

\[\twpa(z,\tau) =  \frac{q^{1/2}}{16\, \sin \left( \pi \,z \right) ^{2}} \prod_{n=0}^{+\infty} (1-q^{2n+2})^4(1-q^{n+1/2})^4 \left( \prod _{n=0}^{+\infty}{\frac { \left( 1-{{\rm e}^{2\,
i\pi \,z}}{q}^{n} \right)  \left( 1-{{\rm e}^{-2\,i\pi \,z}}{q}^{n}
 \right) }{ \left( 1+{{\rm e}^{2\,i\pi \,z}}{q}^{n+1/2} \right) 
 \left( 1+{{\rm e}^{-2\,i\pi \,z}}{q}^{n+1/2} \right) }} \right) ^{2}.\]

\[
\twpa(z,\tau) = \sin \left( \pi \,z \right) ^{2} q^{1/2}
\prod_{n=0}^{+\infty} (1-q^{2n+2})^4 (1-q^{n+1/2})^4
\pa{\ds \prod _{n=0}^{+\infty} \frac{\ds \left( 1-{{\rm e}^{2\,
i\pi \,z}}{q}^{n+1} \right)  \times \\\left( 1-{{\rm e}^{-2\,i\pi \,z}}{q}^{n+1}
 \right) }{\left( 1+{{\rm e}^{2\,i\pi \,z}}{q}^{n+1/2} \right)  \left( 1+{{\rm e}^{-2\,i\pi \,z}}{q}^{n+1/2} \right)}}^2
\]
or, while returning to $\twp_\tau$,
\[
\twp_\tau(z) = -16\pi^2 \sin \left( \pi \,z \right) ^{2} q^{1/2}
\prod_{n=0}^{+\infty} (1-q^{2n+2})^4 (1-q^{n+1/2})^4
\pa{\ds \prod _{n=0}^{+\infty} \frac{\ds \left( 1-{{\rm e}^{2\,
i\pi \,z}}{q}^{n+1} \right)  \left( 1-{{\rm e}^{-2\,i\pi \,z}}{q}^{n+1}
 \right) }{\left( 1+{{\rm e}^{2\,i\pi \,z}}{q}^{n+1/2} \right)  \left( 1+{{\rm e}^{-2\,i\pi \,z}}{q}^{n+1/2} \right)}}^2.
\]


\bs

\textsc{$\bullet$ Factorization of ${\rm SD}(z,\tau)$ :}
\ms

According to the previous section and after standardization, for $\tau\in {\cal H}$ and $z = x + y\tau$:
\[
\begin{array}{ r c l l}
{\rm SD}(z,\tau) & = & \ds \frac{1}{2}\sum_{n=0}^{+\infty} \frac{(-1)^n\sin((2n+1)\pi z)}{\cos((n+\frac{1}{2})\pi\tau)} & \ \ {\rm when} \ \abs{y} < \frac{1}{2}, \\
 & = & \ds -\frac{1}{2\tau}\sum_{n=0}^{+\infty} \frac{(-1)^n \sin((2n+1)\pi z/\tau)}{\cos((n+\frac{1}{2})\pi/\tau)} & \ \ {\rm when} \ \abs{x} < \frac{1}{2}
\end{array}
\]

This is an elliptic function in $z$ (with $\tau$ fixed), $2$ and $2\tau$ periodic. Moreover:
\[{\rm SD}(z+1,\tau) = -{\rm SD}(z,\tau) \ \ {\rm and} \ \ {\rm SD}(z+\tau,\tau) = -{\rm SD}(z,\tau).\]

Then $1+\tau$ and $1-\tau$ are two periods for ${\rm SD}(z,\tau)$.

The function is cancelled only in $\Z+\tau\Z$. This shows that $1+ \tau$ and $1- \tau$ form a fundamental system of periods. In the fundamental domain $[-\frac{1}{2},\frac{1}{2}](1+\tau) + [-\frac{1}{2},\frac{1}{2}](1-\tau)$, $0$ is a zero inside the domain, whereas $-1$, $1$, $\tau$ and $-\tau$ shares each one on four domains.

The relocated domain $\frac{1}{2} + [-\frac{1}{2},\frac{1}{2}](1+\tau) + [-\frac{1}{2},\frac{1}{2}](1-\tau)$ contains exactly $0$ and $1$ as zero which are quite interior.
\ms

Consequently, out of $\Z+\tau\Z$, $SD(z,\tau)$ does not cancel itself. It is the same phenomenon as for $\twpa$, moreover the poles of $SD(z,\tau)^2$ are simple and located on lattice $\Lambda^\star$ and consequently:
\[\forall (z,\tau)\in \C\times{\cal H}, \ \ SD(z,\tau)^2 = \twpa(z,\tau).\]

Indeed, $z\mapsto SD(z,\tau)^2$ is an elliptic function of degree $2$ on the $\Lambda_\tau$ lattice with double-zeros on $\Lambda_\tau$ and poles on $\Lambda_\tau^\star$. It remains to adjust the constant of proportionality with an equivalent in the vicinity of $0$ to conclude with the announced equality using Liouville's theorem.
\ms

Factorization of $\twpa(z,\tau)$ gives the factorization of $SD$ according to $z$ and $q$ variables:

\[SD(z,\tau) =
\sin \left( \pi \,z \right) {q}^{1/4} \prod _{n=0}^{+\infty} \left( 1-{q}^
{2\,n+2} \right) ^{2} \left( 1-{q}^{n+1/2} \right) ^{2}\prod _{n=0}^{N
}{\frac { \left( 1-{{\rm e}^{2\,i\pi \,z}}{q}^{n+1} \right)  \left( 1-
{{\rm e}^{-2\,i\pi \,z}}{q}^{n+1} \right) }{ \left( 1+{{\rm e}^{2\,i
\pi \,z}}{q}^{n+1/2} \right)  \left( 1+{{\rm e}^{-2\,i\pi \,z}}{q}^{n+
1/2} \right) }}.
\]
\bs

\textsc{$\bullet$ The Jacobi functions ${\rm CD}(z,\tau)$ and ${\rm ND}(z,\tau)$ :}
\ms

For recall, when $\tau\in {\cal H}$ with $z = x + y\tau$:

\[
\begin{array}{ r c l l}
{\rm CD}(z,\tau) & = & \ds \frac{i}{2}\sum_{n=0}^{+\infty} \frac{(-1)^n\cos((2n+1)\pi z)}{\sin((n+\frac{1}{2})\pi\tau)} & \ \ {\rm when} \ \abs{y} < \frac{1}{2}, \\
 & = & \ds \frac{i}{4\tau} + \frac{i}{2\tau}\sum_{n=1}^{+\infty} \frac{(-1)^n \cos(2n\pi z/\tau)}{\cos(n\pi/\tau)} & \ \ {\rm when} \ \abs{x} < \frac{1}{2}
\end{array}
\]

\[
\begin{array}{ r c l l}
{\rm ND}(z,\tau) & = & \ds 1 + 2\sum_{n=1}^{+\infty} \frac{(-1)^n\cos(2n\pi z)}{\cos(n\pi\tau)} & \ \ {\rm when} \ \abs{y} < \frac{1}{2}, \\
 & = & \ds \frac{2}{\tau}\sum_{n=0}^{+\infty} \frac{(-1)^n \cos((2n+1)\pi z/\tau)}{\sin((n+\frac{1}{2})\pi/\tau)} & \ \ {\rm when} \ \abs{x} < \frac{1}{2}.
\end{array}
\]

The function $z\mapsto CD(z,\tau)$ is $2$ and $\tau$ periodic, it verifies $CD(z+1,\tau) = - CD(z,\tau)$. Its zero, all simple, are on $\frac{1}{2} + \Z+\tau\Z$. Its poles, all simple, are on $\frac{1+\tau}{2} + 2\Z+\tau\Z$.

The function $z\mapsto ND(z,\tau)$ is $1$ and $2\tau$ periodic and verifies $ND(z+\tau,\tau) = - ND(z,\tau)$. Its zero, all simple, are on $\frac{\tau}{2} + \Z+\tau\Z$. Its poles, all simple, are on $\frac{1+\tau}{2} + 2\Z+\tau\Z$.
\ms


Previous Fourier series representations for $SD$ and $CD$ give the equality:
\[\forall \tau\in {\cal H}, \ \ CD(0,\tau) = -iSD(\frac{1}{2},\tau+1).\]

The zeros and the poles of $CD$ are well located, while following the process of factorization of $ \ twpa$, it comes:

\[CD(z,t) =
\cos \left( \pi \,z \right) {q}^{1/4} \prod _{n=0}^{+\infty}{\frac {
 \left( 1-{q}^{n+1} \right) ^{2}}{ \left( 1-{q}^{n+1/2} \right) ^{2}}}
\prod _{n=0}^{+\infty}{\frac { \left( 1+{{\rm e}^{2\,i\pi \,z}}{q}^{n+1}
 \right)  \left( 1+{{\rm e}^{-2\,i\pi \,z}}{q}^{n+1} \right) }{
 \left( 1+{{\rm e}^{2\,i\pi \,z}}{q}^{n+1/2} \right)  \left( 1+{
{\rm e}^{-2\,i\pi \,z}}{q}^{n+1/2} \right) }}.
\]

From this last equality and factorization of $SD$ one deduces then:
\[\forall \tau\in {\cal H}, \ \ CD(0,\tau) = SD(\frac{1}{2},\tau) \ \equivalent \ i \sum_{n=0}^{+\infty}\frac{(-1)^n}{\sin((n+\frac{1}{2})\pi\tau)} = \sum_{n=0}^{+\infty}\frac{1}{\cos((n+\frac{1}{2})\pi\tau)}.\]

While considering the Fourier series representations, it thus comes:
\[\forall \tau\in {\cal H}, \ \ ND(0,\tau) = 4iSD(\frac{\tau}{2},\tau).\]

\[ND(z,t) =
\prod _{n=0}^{+\infty} \left( 1-{q}^{2\,n+1} \right) ^{2} \left( 1-{q}^{n+1}
 \right) ^{2}\prod _{n=0}^{+\infty}{\frac { \left( 1-{{\rm e}^{2\,i\pi \,z}}
{q}^{n+1/2} \right)  \left( 1-{{\rm e}^{-2\,i\pi \,z}}{q}^{n+1/2}
 \right) }{ \left( 1+{{\rm e}^{2\,i\pi \,z}}{q}^{n+1/2} \right) 
 \left( 1+{{\rm e}^{-2\,i\pi \,z}}{q}^{n+1/2} \right) }}.
\]

Briefly, $\twpa = SD^2$ and one has factorizations:

\begin{equation}
\twpa(z,\tau) = \sin \left( \pi \,z \right) ^{2} q^{1/2}
\prod_{n=0}^{+\infty} (1-q^{2n+2})^4 (1-q^{n+1/2})^4
\pa{\ds \prod _{n=0}^{+\infty} \frac{\ds \left( 1-{{\rm e}^{2\,
i\pi \,z}}{q}^{n+1} \right)  \left( 1-{{\rm e}^{-2\,i\pi \,z}}{q}^{n+1}
 \right) }{\left( 1+{{\rm e}^{2\,i\pi \,z}}{q}^{n+1/2} \right)  \left( 1+{{\rm e}^{-2\,i\pi \,z}}{q}^{n+1/2} \right)}}^2,
\label{fact-twpa}
\end{equation}

\begin{equation}
SD(z,\tau) =
\sin \left( \pi \,z \right) {q}^{1/4} \prod _{n=0}^{+\infty} \left( 1-{q}^
{2\,n+2} \right) ^{2} \left( 1-{q}^{n+1/2} \right) ^{2}\prod _{n=0}^{N
}{\frac { \left( 1-{{\rm e}^{2\,i\pi \,z}}{q}^{n+1} \right)  \left( 1-
{{\rm e}^{-2\,i\pi \,z}}{q}^{n+1} \right) }{ \left( 1+{{\rm e}^{2\,i
\pi \,z}}{q}^{n+1/2} \right)  \left( 1+{{\rm e}^{-2\,i\pi \,z}}{q}^{n+
1/2} \right) }},
\end{equation}

\begin{equation}
CD(z,t) =
\cos \left( \pi \,z \right) {q}^{1/4} \prod _{n=0}^{+\infty}{\frac {
 \left( 1-{q}^{n+1} \right) ^{2}}{ \left( 1-{q}^{n+1/2} \right) ^{2}}}
\prod _{n=0}^{+\infty}{\frac { \left( 1+{{\rm e}^{2\,i\pi \,z}}{q}^{n+1}
 \right)  \left( 1+{{\rm e}^{-2\,i\pi \,z}}{q}^{n+1} \right) }{
 \left( 1+{{\rm e}^{2\,i\pi \,z}}{q}^{n+1/2} \right)  \left( 1+{
{\rm e}^{-2\,i\pi \,z}}{q}^{n+1/2} \right) }},
\end{equation}

\begin{equation}
ND(z,t) =
\prod _{n=0}^{+\infty} \left( 1-{q}^{2\,n+1} \right) ^{2} \left( 1-{q}^{n+1}
 \right) ^{2}\prod _{n=0}^{+\infty}{\frac { \left( 1-{{\rm e}^{2\,i\pi \,z}}
{q}^{n+1/2} \right)  \left( 1-{{\rm e}^{-2\,i\pi \,z}}{q}^{n+1/2}
 \right) }{ \left( 1+{{\rm e}^{2\,i\pi \,z}}{q}^{n+1/2} \right) 
 \left( 1+{{\rm e}^{-2\,i\pi \,z}}{q}^{n+1/2} \right) }}.
\end{equation}
\bs
\ms

{\bf VII - Standard differential equation for $\tilde{\wp}_{\omega_1,\omega_2}$}

\ms

The Weierstrass function $\wp = \wp_{\omega_1,\omega_2}$ is classically the solution of an autonomous differential equation $y'^2 = 4y^3 - 60 g_4 y - 140g_6$, where one posed
$$\forall k\in \N-\{0,1,2\}, \ \ g_k = \sum_{a\in \Lambda_{\omega_1,\omega_2}-\{0\}} \frac{1}{a^k}.$$
\ms

This relation makes it possible in its turn to parameterize elliptic curve by $z\mapsto (\wp(z),\wp'(z))$. The relation $(6)$ shows that $\tilde{\wp}_{\omega_1,\omega_2}$ verifies also a differential equation in the same standard way.

Let us pose $f(z) = \tilde{\wp}_{\omega_1,\omega_2}(z+\omega_0) = \wp_{\omega_1,\omega_2}(z) - \wp_{\omega_1,\omega_2}(\omega_0)$, $f$ is even and has a pole of order $2$ at $0$ with null residue. The formal expansions of $f$ and $f'$ in the vicinity of $0$ are written

$$f(z) = \frac{1}{z^2} + \nu_0 + \nu_2 z^2 + \nu_4 z^4 + O(z^6) \ \ {\rm et} \ \ f'(z) = -\frac{2}{z^3} + 2\nu_2 z + 4\nu_4 z^3 + O(z^5).$$

So, we have $f'^2 - 4 f^3 - af^2 - bf = -28\nu_4 + 20\nu_0 \nu_2 -4\nu_0^3 + O(z^2)$ with $a = -12\nu_0$ and $b = 12\nu_0^2-20\nu_2$. Liouville's theorem makes it possible to assert that $f'^2 - 4 f^3 - af^2 - bf = -28\nu_4 + 20\nu_0 \nu_2 -4\nu_0^3$ and the evaluation at $z = \omega_0$ indicates 
\begin{equation}
-28\nu_4 +20\nu_0 \nu_2 -4\nu_0^3 = 0. \label{eqmod1}
\end{equation}

Thus, $f$ is a solution of the equation
\begin{equation}
y'^2 = 4 y^3 + a(\omega_1,\omega_2)y^2 + b(\omega_1,\omega_2)y. \label{eqdiff1}
\end{equation}

The equation being autonomous, $\tilde{\wp}_{\omega_1,\omega_2}(z)$ is also solution of $\eqref{eqdiff1}$. The relation $\eqref{eqmod1}$ is an equality between modular forms and, in the same way, the nullity of the coefficients of the asymptotic expansion of $f'^2 - 4 f^3 - af^2 - bf$ provides relations between modular forms.
\ms

One notes $\ds \sum {'} = \lim_{N\to+\infty} \sum_{\buildrel {-N\ie n,m\ie N,} \over{(n,m)\not = (0,0)}}$. A direct calculation gives

\[\nu_0 = -\frac{1}{\omega_0^2} + \sum {'} \frac{1}{(n\omega_1+m\omega_2)^2} - \frac{1}{(\omega_0+n\omega_1+m\omega_2)^2} = 4\sum {'} \frac{(-1)^{n+m}}{(n(\omega_1+\omega_2) + m(\omega_1-\omega_2))^2},\]

\[\nu_2 = \ds 3\sum {'} \frac{1}{(n\omega_1+m\omega_2)^4} \ \ {\rm et} \ \ \nu_4 = 5\sum {'} \frac{1}{(n\omega_1+m\omega_2)^6}.\]
Thus,
\[a(\omega_1,\omega_2) = -48 \sum {'} \frac{(-1)^{n+m}}{(n(\omega_1+\omega_2) + m(\omega_1-\omega_2))^2}\]
and
\[b(\omega_1,\omega_2) = 192 \pa{\sum {'} \frac{(-1)^{n+m}}{(n(\omega_1+\omega_2) + m(\omega_1-\omega_2))^2}}^2 - 60 \sum {'} \frac{1}{(n\omega_1+m\omega_2)^4}.\]

It is almost immediate that the coefficients $\nu_{2k}(1,\tau)$ are collinear with the Eisenstein functions $G_{2k}(\tau)$ for $k \se 1$. More surprising, with the representation chosen, $\nu_0(1,2\tau+1)$ is collinear with $G_2^*(\tau)$ [diam], which reinforces the relevance of the convention developed here. Equally, the conjecture abc and the proof by A. Wiles of STW theorem for the semi-stable elliptic curves relate to the cubics such as $y^2 = x^3+\alpha x^2+\beta x$ and not $y^2 = x^3+\alpha x+\beta$ as in the Weierstrass parametrization.
\bs

Finally, following the example of classical case, one can factorize this equation and introduce a discriminant.

The function $\tilde{\wp} = \tilde{\wp}_{\omega_1,\omega_2}$ is even, $\Lambda_{\omega_1,\omega_2}$-periodic, well defined on $\omega_1$ and $\omega_2$, so we have $\tilde{\wp}'(\frac{\omega_1}{2}) = \tilde{\wp}'(\frac{\omega_2}{2}) = 0$.
\ms

On the compact surface $\C/\Lambda_{\omega_1,\omega_2}$ the number of zeros of $\tilde{\wp}'$ is equal to the number of poles (all taking into account the multiplicities). The single pole is of order $3$ located at $\frac{\omega_1+\omega_2}{2}$ modulo $\Lambda_{\omega_1,\omega_2}$, consequently $0$, $\frac{\omega_1}{2}$ and $\frac{\omega_2}{2}$ are the only zero modulo $\Lambda_{\omega_1,\omega_2}$ and they are all of order one.
\ms

And we have :
\[\tilde{\wp}'^2 = 4 \tilde{\wp}^3 + a(\omega_1,\omega_2)\tilde{\wp}^2 + b(\omega_1,\omega_2)\tilde{\wp} = 4\tilde{\wp}(\tilde{\wp} - \tilde{\wp}(\frac{\omega_1}{2}))(\tilde{\wp} - \tilde{\wp}(\frac{\omega_2}{2})).\]

At last, $\tilde{\wp}(z) - \tilde{\wp}(\frac{\omega_1}{2}) = 0$ have exactly two zeros on $\C/\Lambda_{\omega_1,\omega_2}$. Since $\tilde{\wp}'(\frac{\omega_1}{2}) = 0$, we deduce that $\frac{\omega_1}{2}$ is the single zero (of order $2$) of $\tilde{\wp}(z) - \tilde{\wp}(\frac{\omega_1}{2}) = 0$ modulo $\Lambda_{\omega_1,\omega_2}$.

Accordingly $\tilde{\wp}(\frac{\omega_1}{2})\not= \tilde{\wp}(\frac{\omega_2}{2})$.

The discriminant of the equation $\tilde{\Delta}(\omega_1,\omega_2) = a^2 - 16b = 16(\tilde{\wp}(\frac{\omega_1}{2})-\tilde{\wp}(\frac{\omega_2}{2}))^2$  is of degree $4$ and never cancels itself.
\bs

When we define the reduced discriminant

\[\forall \tau\in {\cal H}, \ \ \tilde{\Delta}(\tau) = \tilde{\Delta}(1,\tau),\]

The equalities for $\tau \in {\cal H}$:

\[\tilde{\Delta}(\tau+2) = \tilde{\Delta}(\tau) \ {\rm et} \ \tilde{\Delta}(-\frac{1}{\tau}) = \tau^4\tilde{\Delta}(\tau)\]

allow to pose $\Delta_4(\tau) = \tilde{\Delta}(2\tau)$.

By usual techniques, the previous equalities show that $\Delta_4$ is a modular form of weight $4$ and level $4$, which is not cancelled on $\cal H$. It plays on its space the same usual part as discrimant $\Delta$ on the space of the modular forms of weight $12$ and level $1$.
\ms \bs

{\bf VIII - Inversion of $\tilde{\wp}_{\omega_1,\omega_2}$}
\ms

By construction, $\tilde{\wp}_{\omega_1,\omega_2}$ is an even function, it presents one double-zero at $0$ and one double-pole at $\omega_0 = \frac{\omega_1+\omega_2}{2}$, with a null residue.

Thanks to Liouville's theorem, the product $z\mapsto \tilde{\wp}_{\omega_1,\omega_2}(z)\tilde{\wp}_{\omega_1,\omega_2}(z-\omega_0)$ is a non-null constant function. Let us note $C(\omega_1,\omega_2)$ this value.

Using parity, we have
\[C(\omega_1,\omega_2) = \tilde{\wp}_{\omega_1,\omega_2}(\frac{1}{4}(\omega_1+\omega_2))^2 = \tilde{\wp}_{\omega_1,\omega_2}(-\frac{1}{4}(\omega_1+\omega_2))^2.\]
\ms

The equation \eqref{eqdiff1} indicates that $\ds \frac{1}{\tilde{\wp}_{\omega_1,\omega_2}} = \ds \frac{f}{C}$ is solution of $y'^2 = b(\omega_1,\omega_2)y^3 + a(\omega_1,\omega_2)y^2 + 4y$. According to what precedes, the function $f$ being solution only of one single equation of this type, i.e. \eqref{eqdiff1} itself, it thus comes

\[C(\omega_1,\omega_2) = \frac{1}{4}b(\omega_1,\omega_2).\]
\ms

The function $C(\tau) = C(1,\tau)$ is $2$ periodic and verifies $C(\tau) = \tau^4C(-1/\tau)$. One poses finally $C_1(\tau) = C(2\tau)$ and, following the example of $\Delta_4$, one shows that $C_1$ is a modular form of weight $4$ and level $4$ which is not cancelled on $\cal H$.


\ms\bs



{\bf IX -- Conclusion and prospects}
\ms

The study of the elliptic functions, according to the approach of Weiertrass and Jacobi, forms a mathematical theory whose basic elements are standardized. For example Weierstrass's $\wp$ function presents a double-pole at the origin. Using contour integrals, whose introduction is justified by a global approach of the summation formulae [Feau], we found rather simply almost the usual elliptic functions.
\ms

This translation appears naturally in the calculation of the contour integrals, and this missing element enabled us to factorize the renormalized elliptic functions.
\ms

This renormalization led to Fourier series representations of the Weierstrass or Jacobi elliptic functions on the reduced lattice, and the observation immediately raised a modular character.
\ms

These Fourier series expansions will be an access point to certain modular functions. They will allow, for example, to propose bases of modular spaces whose elements are not cancelled on the half-plane of Poincaré (in fact, modular units). Some examples charateristic of these objects appeared besides naturally in the study of the differential equation related to a Weierstrass's type elliptic function, or in the inversion of such a function.
\ms

Regarding the applications of the formulae of factorization, some are straightforward, here is an example.

One points out the notations $\ds \twpa(z,\tau) = -\frac{1}{16\pi^2}\twp_\tau(z)$ and $q = e^{2i\pi\tau}$, a direct application of the formula of factorization \eqref{fact-twpa} leads successively to:

\[-\frac{1}{16\pi^2}\twp_\tau(\frac{1}{2}) =
q^{1/2}\prod _{k=0}^{+\infty}{\frac { \left( 1-{q}^{2\,k+2} \right) ^{4}}{
 \left( 1-{q}^{2\,k+1} \right) ^{4}}}
,\]

\[\frac{1}{\pi^2}\twp_{t}(\frac{\tau}{2}) =
\prod _{k=0}^{+\infty} \left( 1-{q}^{k+1/2} \right) ^{8} \left( 1-{q}^{k+1}
 \right) ^{4}
,\]

\[\frac{1}{\pi^2}\twp_{\tau+1}(\frac{\tau+1}{2}) =
\prod _{k=0}^{+\infty}{\frac { \left( 1-{q}^{2\,k+1} \right) ^{12} \left( 1-
{q}^{2\,k+2} \right) ^{4}}{ \left( 1-{q}^{k+1/2} \right) ^{8}}}
.\]

Consequently, the modular discriminant $\ds \Delta(q) = q\prod_{k=1}^{+\infty}(1-q^k)^{24}$ can be factorized:

\[\Delta(q) = \frac{1}{256\pi^{12}} \pa{\twp_\tau\pa{\frac{1}{2}}\twp_\tau\pa{\frac{\tau}{2}}\twp_{\tau+1}\pa{\frac{\tau+1}{2}}}^2.\]
\bs
\ms

The author expresses his gratitude to Annie Slark for reviewing the English version of the article.
\bs
\bs

{\textsc{References :}}\ms

\bib{Halp}{G.H. Halphen, {\it Traité des fonctions elliptiques et de leurs applications}, Jacques Gabay, 1886.}

\bib{Armi}{V. Armitage et W. F. Eberlein, {\it Elliptic Functions}, Cambridge University Press, 2006.}

\bib{Apos}{T. M. Apostol, {\it Modular Functions and Dirichlet Series in Number Theory}, Springer-Verlag New York, 1976.}

\bib{Whit}{E.T Whittaker, G.N Watson, {\it A Course of Modern Analysis}, Cambridge University Press, London, 1940.}

\bib{Diam}{F. Diamond, J. Shurman, {\it A First Course in Modular Forms}, Springer-Verlag New York, 2005.}

\bib{Feau}{J.-C. Feauveau, {\it A unified approach for summation formulae}, eprint arXiv:1604.05578}


\bib{Rous}{J. Rouse, J. J. Webb, {\it On spaces of modular forms spanned by eta-quotients}, eprint arXiv:1311.1460}

\bib{Kohl}{G. Köhler, {\it Eta Products and Theta Series Identities}, Springer-Verlag Berlin, 2011.}

\end{document}